\documentstyle[12pt,sr-stylefile]{article}
\renewcommand{\thefootnote}{}

\newcommand{\E}{{\Bbb E}}

\newcommand{\RR}{{ \underline{\cal R}   }}
\newcommand{\RRR}{ {\cal R}  }
\newcommand{\half}{{  {1\over 2}  }}

\newcommand{\ssp}{ {\rm s.~s.~p.}}
\newcommand{\sqq}{ \sqrt{2} }
\newcommand{\osq}{\over \sqrt{2} }
\newcommand{\ft}{\frac{4}{3}}

\def\Bbb{\bf}
\def\limsup{\mathop{\overline{\rm lim}}}
\def\liminf{\mathop{\underline{\rm lim}}}
\def\exp{{\rm e}}

\begin{document}



\Large
\centerline{\bf Bounded and $L^2$ Harmonic Forms on Universal Covers}
\bigskip\bigskip
\large
\centerline{\bf K. D. Elworthy
}
\normalsize 
\medskip
\centerline{\parbox{5in}{Mathematics Institute, University of Warwick,
Coventry CV4 7AL, UK, kde@maths.warwick.ac.uk}}

\bigskip
\large
\centerline{\bf Xue-Mei Li}
\normalsize
\medskip
\centerline{\parbox{5in}{Department of Mathematics,
University of Connecticut,
Storrs, CT 06269, USA, xmli@math.uconn.edu}}

\bigskip
\large
\centerline{\bf Steven Rosenberg}
\normalsize
\medskip
\centerline{\parbox{5in}{
 Department of Mathematics, Boston University,
Boston, MA 02215, USA, sr@math.bu.edu}}

\bigskip
\bigskip





\footnote{The first two authors
are supported by SERC grant GR/H67263.  The last two authors are partially
supported by the NSF.}

\section{Introduction}

In this paper we relate certain curvature conditions on a complete
Riemannian manifold to the existence of bounded and $L^2$ harmonic
forms.  In the case where the manifold is the universal cover of a
compact manifold, we obtain topological and geometric information
about the compact manifold.  Many of these results can be thought of
as differential form analogues of Myers' theorem that a complete manifold
with strictly positive Ricci curvature is compact with finite
fundamental group.   In fact, 
our curvature conditions involve the positivity of
the Weitzenb\"ock curvature terms $\RRR^p$, where $\RRR^1$ is the Ricci
curvature.
We also give pinching conditions on certain sums
of sectional curvatures which imply the positivity of $\RRR^p$, and
hence yield
Bochner-type 
vanishing theorems for harmonic
forms.  In particular, we construct a compact manifold with
planes of negative sectional curvature at each point and which
satisfies the hypothesis of our vanishing theorem.  
Finally, we remark on
related spectral gap 
estimates for the Laplacian on forms on complete
manifolds.  

In more detail, recall that 
Bochner's theorem states that there are no $L^2$ harmonic $p$-forms on
a
complete manifold $M$
if the curvature term $\RRR^p$ in the
Weitzenb\"ock decomposition $\Delta^p = \nabla^\ast\nabla + \RRR^p$
for the Laplacian on $p$-forms is positive.  In particular, if $M$ is
compact, then $\RRR^p >0 \Rightarrow H^p(M;{\bf R}) =0.$  
This result has the
disadvantage that $\RRR^p$ is a complicated curvature expression whose
geometric/topological significance is unclear, 
except for $\RRR^1 = {\rm
Ric}$, the Ricci curvature.  In general, even manifolds with positive
sectional curvature need not have $\RRR^p >0$ for $1<p<n-1.$

We have previously made some progress in this area.  In \cite{er-inv},
we showed that $\RRR^p$ need only be strongly stochastically positive,
or s.~s.~p.,
(as defined in \S2) for Bochner's theorem to hold.  In practice, this
allows $\RRR^p$ to be negative on a set of small volume if $M$ is
compact.
We also
showed that $\RRR^2 $ s.~s.~p.  implies
$\pi_2(M)$ is a torsion group.  In \cite{EL-RO94}, we showed that the
positivity of $\RRR^p$ on primitive vectors is equivalent to the
positivity of a sum of sectional curvatures, and we used this to
obtain vanishing theorems for $H_p(M;{\bf Z})$ for manifolds
isometrically immersed in ${\bf R}^N.$  Finally, in \cite{group} the
second author showed that the hypothesis of Myers' theorem can be
extended to the
strong stochastic positivity of the Ricci
curvature.  

In this paper, we extend these results in several directions.  First,
in contrast to the usual Jacobi field proof,
the proof of Myers' theorem in \cite{group} comes down to showing the
existence of an $L^2$ harmonic function on the universal cover of
$M$.  In \S2, we extend this argument from functions to forms.  The
main analytic result (Theorem 2.1) states that if $\RRR^{p\pm 1}$ is
$\ssp$, and if $M$ has a bounded harmonic
$p$-form, then $M$ has an $L^2$ harmonic $p$-form. 
In fact, all the results in this section carry over to the
Bismut-Witten Laplacian, and are stated in this generality.

As consequences, we give a series of results indicating the 
geometric significance of $\RRR^2.$  
We show that  the universal cover of a 
compact manifold $M$ with  $\RRR^2$ $\ssp$ and with a nonconstant
bounded harmonic function 
admits a nonconstant harmonic function of bounded
energy (Proposition 2.1). This implies that a compact manifold
cannot admit both a
metric with $\RRR^2$ $\ssp$ and a metric with pinched negative
curvature (Theorem 2.3).  (There are corresponding results for
$p$-forms.)  Moreover, if $M,N$ are compact manifolds with $\pi_1(M)$
nonamenable and $\pi_1(N)$ infinite, then $M\times N$ admits no metric
with $\RRR^2 \ \ssp$ (Theorem 2.4).
We also show
that a 4-manifold with nonamenable fundamental group and $\chi(M) >-2$
has no metric with $\RRR^2$ $\ssp$  (Proposition 2.2); in particular,
no $4$-manifold of negative curvature admits such a metric.

In \S3, we prove that $\RRR^p >0$ if there exists $A = A(x)>0$ such that
\begin{equation}\label{first}
C  A < \sum_{i=1}^p \sum_{j=p+1}^n K(v^i, v^j)\equiv \mathop{\sum_p}
 < A,\end{equation}
for all orthonormal bases $\{v^1,\ldots, v^n\}$ of $T_x M$.  Here
$C= C(n,p)$ is an explicit constant, and
$K(v^i, v^j)$ is the sectional curvature of the $(v^i, v^j)$-plane.
In particular, a manifold with this pinching estimate has no $L^2$
harmonic $p$-forms.

Note that the sum of curvatures in (\ref{first}) 
is ${\rm Ric}( v^1,
v^1)$ if $p=1.$  In general, this sum of curvatures is precisely
the sum in \cite{EL-RO94} mentioned above (which in turn is based on 
\cite{lawson-simons}).  Combined with Theorem 2.3,
this shows that a compact manifold cannot admit both a metric of
pinched negative curvature and a metric with $\mathop{\sum_2}$ pinched
as in (\ref{first}).
We also show that the
metric product $M = \Sigma \times S^4$ of a surface $\Sigma$ of 
constant negative curvature with the 4-sphere satisfies the pinching
estimate above with $p=3.$  As mentioned above, this example is
significant because it has planes of negative sectional curvature at
each point.

Finally, in \S4 we give  estimates for the spectral gap at zero for
the (Bismut-Witten) Laplacians on forms.  These estimates come from the $L^2$
analogues of the $L^\infty$ estimates of \S2.

We would like to thank Feng-Yu Wang for helpful conversations.

\section{Bounded harmonic forms and $L^2$ harmonic forms}

In this section we relate the positivity of the Weitzenb\"ock term for
the Laplacian on forms to the existence of bounded and $L^2$ harmonic
forms.  We give applications of this result to the geometry and
topology of compact manifolds.

Let $(M,g)$ be a complete Riemannian manifold and $h:M\to {\bf R}$
 a smooth function. The case $h\equiv 0$ is of particular
 interest. Let $\delta^h$ be the adjoint of the exterior derivative
 $d$ on forms with respect to the measure $\mu(dx)=\exp^{2h}{\rm
 dvol}_g$. The
Bismut-Witten Laplacian $\Delta^h=(d+\delta^h)^2$ is a
self-adjoint operator on $L^2$ forms with respect to
$\mu(dx)$, see \cite{Li.thesis}.  This Laplacian restricts to an
 operator $\Delta^{h,q}$ on $q$-forms, and has the Weitzenb\"ock
 decomposition 
$$\Delta^{h,q}= \nabla^\ast\nabla+{\cal R}^{h,q}.$$
In particular, $\RRR^{h,1} = {\rm Ric} - 2\cdot {\rm Hess}(h)$, where
${\rm
Hess}(h) = -\nabla^2 h$ is the Hessian of $h$.
We use
the convention $\RRR^{h,-1}=0$, $\RRR^{h,n+1}=0$.
Let
$\RR^{h,q}(x)$ be the infimum of ${\cal R}^{h,q}(v)$ over all unit
 $q$-covectors $v\in \Lambda^q T_x^\ast M.$
We will omit the superscripts $h,q$ if the context is clear.

Let $\{x_t\}$ be a path continuous diffusion process with $-\Delta^h =
-\Delta^{h,0} $
as
generator, i.e. a $h$-Brownian motion, starting from $x_0$. Assume
that the 
$h$-Brownian motion does not explode, or equivalently that the
associated Bismut-Witten heat semigroup
 $P^h_t$ is conservative: $P^h_t1\equiv 1$. 
For $P_t^{h,q}=\exp^{-\half t\Delta^{h,q}}$ the heat semigroup on
$L^2\Lambda^q(M, \mu(dx))$ defined by the spectral theorem, 
it is known that
 for each $q$-form $\phi\in L^\infty\cap L^2$,
\begin{equation}\label{fk}
P_t^{h,q}\phi(v_0)=\E \phi\left(W_t^{h,q}(v_0)\right),\end{equation}
provided   $\RR^{h,q}$ is bounded from below. Here $W_t^{h,q}(v_0)$
 is the solution to the stochastic covariant differential equation
 along paths of $\{x_t\}$:
\begin{equation}\label{f-k}
\left\{ \begin{array}{lll}
{DW_t^{h,q}\over \partial t}&=&-\half
{\cal R}^{h,q}( W_t^{h,q}(v_0))\\
W_0^{h,q}(v_0)&=&v_0,
\end{array}\right.\end{equation}
for $v_0\in \Lambda^q T_{x_0}M$; we use the metric $g$ to 
let $\RRR^{h,q}$ act on $\Lambda^q T_{x_0}M.$
 This easily implies
\begin{equation}\label{Witten1}
|W_t^{h,q}(v_0)|\le \exp^{-\half \int_0^t \RR^{h,q}(x_s)ds}.
\end{equation}

\noindent   Let
$$\RR_q(x_0)
=\int_0^\infty \E \exp^{-\half \int_0^t \RR^{h,q}(x_s)ds}dt, $$
and let
$H= H^{h,q}$ be the $h$-harmonic projection on the space of
$L^2\Lambda^q(M,\mu(dx))$.
Thus for such
$\phi$, $H\phi=\lim_{t\to \infty}P_t^{h,q}\phi$.  Let
$C_c^\infty\Lambda^q = C_c^\infty\Lambda^q T^\ast M$ denote the space
of smooth compactly supported $q$-forms, and let $|\phi|_{L^p}$
denoted the $L^p$ norm of $|\phi|.$

\begin{lemma} \label{le:1}
Let $q\in\{2,\ldots,n-2\}.$
 Assume that the $h$-Brownian motion does not explode and that
 ${\cal R}^{h,q-1}$ and  ${\cal R}^{h,q+1}$ are bounded below.
For  $\phi, \psi\in C_c^\infty\Lambda^q$, we have
\begin{eqnarray}
\left|\int_M\langle H\phi-\phi,\psi\rangle\mu(dx)\right|
&\le&
 \half\left[\sup_{x_0\in
 Supp(\phi)}\RR_{q+1}(x_0)\right]|d\psi|_{\infty}
\cdot |d\phi|_{L^1}\\
&&\qquad + \half \left[\sup_{x_0\in Supp(\phi)}\RR_{q-1}(x_0)\right]
|\delta^h\psi|_{\infty}\cdot |\delta^h\phi|_{L^1}.\nonumber
\end{eqnarray}
\end{lemma}

\noindent{\sc Proof.}
For $P_t = P_t^{h,q}$, we have
\begin{equation}\label{10}
\int_M\langle H\phi-\phi,\psi\rangle  \mu(dx)=\lim_{t\to \infty}
\int_M\langle P_t\phi-\phi,\psi\rangle \mu(dx).
\end{equation}
On the other hand,
\begin{eqnarray*}
\int_M\langle P_t\phi-\phi,\psi\rangle \mu(dx)
&=&-\half\int_M \int_0^t \langle \Delta^h(P_s\phi),\psi\rangle
\mu(dx)\\
&=&\half \int_0^t \int_M \langle d\phi, P_s(d\psi)\rangle \mu(dx)
\,ds\\
&&\qquad
+\half \int_0^t \int_M \langle \delta^h\phi, P_s(\delta^h \psi)\rangle
\mu(dx)\, ds.
\end{eqnarray*}

\noindent Under the assumptions of the lemma,
 $|P_s(d\psi)|(x_0)\le |d\psi|_\infty\cdot \E|W_s^{h,q+1}|_{x_0}$
and 
$|P_s(\delta^h(\psi)|(x_0)\le |\delta \psi|_\infty\cdot
\E|W_s^{h,q-1}|_{x_0}$.
So by (\ref{Witten1}) we have
\begin{eqnarray*}
\lefteqn{\left|\int_M\langle P_t\phi-\phi,\psi\rangle
\mu(dx)\right|}\\
&&\le
 \half |d\psi|_\infty \cdot \sup_{x_0\in supp(\phi)} \int_0^\infty 
\E \exp^{-\half \int_0^u \RR^{h,q+1}(x_s)ds}du\cdot|d\phi|_{L^1}\\
&&\qquad+\half |\delta^h\psi|_\infty \cdot
 \sup_{x_0\in {supp}(\phi)}\int_0^\infty 
 \E \exp^{-\half \int_0^u \RR^{h,q-1}(x_s)ds}du
\cdot|\delta^h\phi|_{L^1},
\end{eqnarray*}
and the required inequality follows from (\ref{10}).
 \hfill\rule{3mm}{3mm}

\bigskip

 For one-forms the corresponding result is:

\begin{lemma}\label{le:2}
Assume that the $h$-Brownian motion does not explode. Then
for $\phi, \psi\in  C_c^\infty\Lambda^1T^\ast M$,
\begin{eqnarray*}
\int_M\langle H\phi-\phi,\psi\rangle \mu(dx)
&\le&
\half \left[\sup_{x_0\in Supp(\phi)}\RR_{2}(x_0)\right] 
|d\psi|_\infty \cdot|d\phi|_{L^1}\\
&&\qquad +\half \sup_{x_0\in {supp}(\phi)} \biggl|\int_0^\infty
P_s(\delta^h\phi)ds\biggr|\cdot |\delta^h \psi|_{L^1}
\end{eqnarray*}
\end{lemma}

The main application of these estimates 
is to the existence of harmonic forms. By a 
{\it   bounded $C^1$ h-harmonic  $q$-form} ($q\neq 0$), we mean a
bounded $q$-form  
with $d\phi=0$ and $\delta^h\phi=0$.  For $q=0$, a bounded harmonic
function $f$ is given by the usual definition $\Delta^{h,0} f =0.$  
Note that a bounded solution of
$\Delta^{h,q}\phi =0$ need not be a bounded $h$-harmonic form.  By an
$L^2$ {\it h-harmonic q-form} we mean an $L^2$ solution of
$\Delta^{h,q}\phi =0$.  An $L^2$ $h$-harmonic $q$-form $\phi$
satisfies $d\phi = \delta^h\phi =0.$  The following result is based on
techniques going back to \cite{bakry}.

\begin{theorem}\label{pr:1}
Let $q\in \{ 2,3,\dots, n-2\}$. Assume that h-Brownian motion does not
 explode.  Suppose ${\cal R}^{h,q-1}$ and ${\cal R}^{h,q+1}$ are
 bounded from below with both $\sup_{x\in K}\RR_{q+1}(x)$ and
 $\sup_{x\in K}\RR_{q-1}(x)$ finite for all compact
sets $K\subset M.$  If there exists a nontrivial
 bounded h-harmonic $q$-form, then there exists a nontrivial $L^2$
 harmonic
 $q$-form.
\end{theorem}

\noindent{\sc Proof.}
Let $\psi$ be a nontrivial bounded $h$-harmonic $q$-form.  
Take a sequence of
functions $\{h_n\}$ on $M$ with $0\le h_n \le1$,
 $|\nabla h_n|\le {1\over n}$ and $h_n\to 1$ as $n\to \infty$,  
see \cite{bakry}, \cite{gaffney}.  Define $\psi_n=h_n\cdot \psi$. Then
$d\psi_n=d h_n \wedge \psi$ and $\delta^h \psi_n=-{\rm i}_{\nabla h_n}\psi$,
and so 
$$|d\psi_n|_\infty\le {1\over n}|\psi|_\infty, \hskip 12pt
|\delta^h \psi_n|_\infty\le {1\over n} |\psi|_\infty.$$
Thus by Lemma \ref{le:1},
$$\left|\int_M \langle H\phi-\phi, \psi_n\rangle  \mu(dx) \right|\le
 c{|\psi|_\infty\over n}
 \left[|d\phi|_{L^1}+|\delta^h\phi|_{L^1}\right],$$
where $c$ is a constant depending on the support of $\phi$. If there
 are no $L^2$ $h$-harmonic $q$-forms, then
$H\phi\equiv 0$ for all $\phi\in C_c^\infty\Lambda^q$. Thus
$$\lim_{n\to \infty}\int_M\langle -\phi, \psi_n\rangle \mu(dx)
=-\int_M \langle \phi,\psi\rangle \mu(dx),$$
since $\psi$ is bounded and $\phi$ has compact support. So
$$\int_M \langle \phi,\psi\rangle \mu(dx)=0$$
for any $\phi\in C_c^\infty\Lambda^q$. Therefore $\psi$ has to be
zero, 
contradicting the hypothesis.\break
${}$\hfill\rule{3mm}{3mm}
\bigskip

\noindent{\bf Definition:} A function $f\in C^1(M)$ is {\it strongly
stochastically positive (s.~s.~p.)} if
$$\limsup_{t \to \infty}{1\over t} \sup_{x_0\in K} \log
\E\left(\exp^{-\half \int_0^t f(x_s)ds}\right)<0$$
for each compact set $K$.  We say that $\RRR^p$ is $\ssp$ if the
function $\underline{\RRR}^p$ is $\ssp.$

\bigskip
\noindent {\bf Remarks:}  (1) Strong stochastic positivity is
equivalent to $\Delta^0 + f >0$ on
$L^2(M,g)  $ if $M$ is compact 
\cite{EL-RO88}.  In particular, a positive function on a compact
manifold is
strongly stochastically positive.  More generally, if $f$ is positive
except on a set of small volume in the sense of \cite{er-inv},
\cite{ry}, then $f$ is $\ssp$ Note that a function is $\ssp$ on any
manifold $M$, then its pullback is $\ssp$ on any cover of $M$, since
Brownian motion on the cover projects to Brownian motion on the base.
Thus if  $\Delta^0 + 
\underline{\RRR}^p >0$ on
  a compact manifold $M$,
then $\RRR^p$ is $\ssp$ on the universal cover
$\widetilde M$ with respect to the pullback metric.  Thus the
condition $\RRR^p$ is $\ssp$ on $\widetilde M$ generalizes the
condition $\RRR^p >0$ on $M$.

 (2) By semigroup domination \cite{EL-RO88},
$\RRR^{h,q}$ $\ssp$ implies the vanishing of the spaces of $L^2$
$h$-harmonic $q$-forms and bounded $h$-harmonic $q$-forms.
  Note also that $\RRR^{h,q}$ $\ssp$ implies
that $\sup_{x\in K} \RR_q(x)$ is finite, for $K$ compact in $M$.
Thus the hypotheses in
Theorem 2.1 on this quantity can be replaced by assuming $\RRR^{h,q\pm
1}$ is  $\ssp$  

\bigskip
In summary, {\it all hypotheses in this paper of the form ``$\RRR^p$ is
$\ssp$ on the compact manifold $M$'' can be replaced by ``$\RRR^p >0$ on
$M$.''}

\bigskip

We now give some geometric applications of the last theorem. 
Let  $H^k(M)$
and  $ L^2H^k(M)$ denote the  spaces of  harmonic 
   and $L^2$  $h$-harmonic $k$-forms (or harmonic $k$-forms in the
case of  
  $h\equiv 0$).  Of course, if $M$ is compact then $H^k(M) \simeq
  H^k(M;{\bf R}).$

 A  manifold is
said to  have {\it  negative curvature with pinching constant $k$}
if  its sectional  curvatures $K$ satisfy   $-1\le K<-k<0.$

\begin{corollary}\label{co:pincing1}
Let $M$ be a compact n-manifold.  Suppose $\RRR^{q\pm 1}\equiv
 \RRR^{0,q\pm 1}$ are strongly stochastic positive for some q with
 $q\in \{2, \ldots,n-2\}$ and $q\neq n/2.$
 If $H^q(M)\not =0$, then $M$ admits no metric of
 negative curvature with pinching constant $q^2/(n-q-1)^2$.
\end{corollary}

\noindent {\sc Proof.}
As noted in Remark (1), 
strong stochastic positivity is preserved in passing to a
cover.  In particular,
the universal cover   $\widetilde M$ of
$M$ has $\RRR^{q\pm 1}$ $\ssp$  Since
 $H^q(M;{\bf R})\neq 0$, there is a nontrivial harmonic $q$-form on
$M$, which lifts to a bounded harmonic form on $\widetilde M.$
By Theorem 2.1,
$L^2H^q(\widetilde M) \neq 0.$

Suppose  $M$ also admits a metric with
negative curvature pinched as in the hypothesis.
 Then  $L^2H^q(\widetilde M)=0$, with respect to
this new metric \cite{EL-RO93}. 
However the vanishing of $L^2H^q(\widetilde M)$ is independent of the
 Riemannian metric on $M$ \cite{atiyah}, which gives a contradiction. 
\hfill\rule{1.5mm}{3mm}

\bigskip

\noindent {\bf Remark:} Theorem 2.1 generalizes the corresponding
 result for $q=0$ in \cite{group}.  In particular, since every
 manifold
 has a bounded harmonic function, namely $1$,  the universal cover
 $\widetilde M$ of a compact manifold $M$ with strongly stochastic
 positive Ricci curvature has an $L^2$ harmonic function.  Such
 functions are constant, which
 implies
 that $\widetilde M$ has finite volume. Thus 
 $M$ must be finite volume and have finite
 fundamental group.  This extends Myers' theorem, which states
 that a complete manifold with ${\rm Ric} >\epsilon >0$ is compact
 with finite fundamental group.
\bigskip

Next we discuss the corresponding result for differential
one-forms. Recall
that our diffusion process is transient if
$$\sup_{x\in K} \int_0^\infty P_tf(x)dt<\infty$$
for each $f$ with compact support, see e.g. \cite{azencott}.

\begin{theorem}\label{one-forms}
Assume the h-Brownian motion does not explode and is transient.
Suppose 
$\RR^{h,2}$ is bounded from below and $\sup_{x\in K} \RR_2
(x)<\infty$ for each
 compact subset $K$. 
 If there exists a nontrivial bounded h-harmonic one-form, then
 there exists a nontrivial $L^2$ h-harmonic one-form.
\end{theorem}


The proof is similar to that of Theorem 2.1.  This leads to the
following stronger result:



%

\begin{proposition} Let $M$ be a compact manifold  
admitting a metric with $\RRR^2$ strongly stochastic
positive.  If either

(i) $\widetilde M$ admits a nonconstant bounded harmonic function with
respect to the pullback metric, or

(ii)
$H^1(M;{\bf R}) \neq 0$ and Brownian motion is transient on
$\widetilde M$,

\noindent then
$L^2H^1(\widetilde{M})\not= 0$, 
and $\widetilde M$ admits a harmonic function of finite energy.
\end{proposition}
\begin{theorem} A compact $n$-manifold cannot admit a metric of 
$[1-(n-2)^{-2}]$-pinched
negative curvature and a metric with $\RRR^2$ strongly stochastically
positive.  \label{cannot}
\end{theorem}

\medskip
\noindent{\sc Proof of the Theorem:}
By \cite{EL-RO93}, the pinching condition on the negatively curved
metric $g_1$ implies 
$L^2H^1(\widetilde{M}, \tilde g_1) =0$, where $\tilde g_1$ is the
pullback of $g_1$ to $\widetilde M.$   A negatively curved manifold
$M$ has
$\lambda_0(\widetilde{M}) >0$ by Mc\-Kean's estimate \cite{mckean},
 where $\lambda_0$
is the bottom of the spectrum for $\Delta^0$. By Brooks' theorem
\cite{brooks}, 
$\pi_1(M)$ must be nonamenable. By \cite{ls}, the nonamenability of
$\pi_1(M)$ implies
$\widetilde{M}$ has a nonconstant bounded harmonic function $f$, which
can be assumed to be positive, with respect to any pullback metric, in
this case the pullback of the metric $g_2$ with $\RRR^2$ $\ssp$ The
the condition
$L^2H^1(\widetilde{M})=0$
is independent of the metric on $M$, since the spaces of $L^2$
harmonic forms are isomorphic to the $L^2$ closed forms modulo the
closure of the $L^2$ exact forms, and hence are quasi-isometry
invariants of $\widetilde M.$  Thus
$L^2{
H}^1(\widetilde{M},\tilde g_2)=0$. This contradicts the
proposition. \hfill\rule{1.5mm}{3mm}

\medskip
{\sc Proof of the Proposition:} For (i), it
 is well known that the existence of a
nontrivial bounded harmonic function $f$ implies 
the transience of  Brownian motion on $\widetilde M$ with respect to
the pullback metric.
 Thus $\omega =df$ is a nonzero one-form
satisfying $d\omega =\delta \omega =0$. Moreover, $\omega$
is bounded by Harnack's inequality \cite{yau}
 or by Bismut's formula in
\cite{el-li94}.  Thus Theorem \ref{one-forms} implies $L^2{
H}^1(\widetilde{M})\not= 0$.  Moreover, $f$ has $\Delta f =0$ and has
 finite energy $\Vert df\Vert_{L^2}.$

For (ii), note that $H^1(M;{\bf R}) \neq 0$ implies
that $M$ has a nontrivial harmonic one-form, which pulls back to a
bounded harmonic one-form on $\widetilde M.$  By Theorem 2.2,
$\widetilde M$ has an $L^2$ harmonic one-form $\omega.$ Since
$H^1(\widetilde M;{\bf R}) =0$, we have $\omega = df$ for some
function $f$ on $\widetilde M.$  Thus $0 = \delta \omega = \Delta f$,
and $f$ has finite energy as before.
\hfill\rule{1.5mm}{3mm}
\bigskip

As an example of the sharp
 nature
of this result,
let $M=S^3\times S^1$. Since the standard metric on 
$S^3$ has $\RRR^2 >0$, the product of the standard
metrics has $\RRR^2>0$ by an easy calculation.
$M$ has a harmonic one-form, which lifts to a
bounded harmonic one-form on $\widetilde{M} = S^3 \times {\bf R}$.
However, $L^2{
H}^1(\widetilde{M})=0$ by direct calculation. Brownian motion on
$\widetilde{M}$ is recurrent, so the hypothesis on the transience of
Brownian motion cannot be omitted from Theorem \ref{one-forms}.

As another application of Theorem \ref{one-forms}, we have

\begin{theorem}  Let $M,N$ be compact manifolds with $\pi_1(M)$ 
nonamenable and $\pi_1(N)$
 infinite.  Then $M\times N$ admits no metric with $\RRR^2 $ strongly
stochastically positive.  \label{product} \end{theorem}

\noindent{\sc Proof:}  The universal covers $\widetilde M, \widetilde N$ have 
no $L^2$ harmonic functions since $\pi_1(M), \pi_1(N)$ are infinite.  The 
K\"unneth formula for $L^2$ harmonic forms thus implies $L^2H^1(\widetilde
 M\times \widetilde N) =0.$  However, $\pi_1(M\times N)$ is nonamenable, so as
 above $\widetilde M\times \widetilde N$ has a nontrivial bounded harmonic 
one-form with respect to any pullback metric.  The theorem follows as above.
\hfill\rule{1.5mm}{3mm}

\bigskip

\noindent{\bf Remark:} As pointed out to us by J. Lott,
the K\"unneth formula for the Laplacian on $L^2$ forms on
a product manifold
follows from the discussion of the spectrum of
an operator of the form $A\otimes 1 + 1\otimes B$ in
\cite[Vol. I, Thm. VIII.33]{reed-simon}.  Namely, this theorem implies that the
zero spectrum of the Laplacian on a product manifold is given by the expected
combination of the zero spectra on the individual factors.

\bigskip

We now consider these results for $4$-manifolds.
In this dimension, the positivity (resp. $\ssp$)
 of $\RRR^2$ is equivalent to the 
positivity (resp. $\ssp$)
 of the curvature operator
on complex isotropic two-planes \cite{m-w}.
The Euler characteristic of such
manifolds is constrained as
follows:

\begin{proposition}\label{fourdims}
Let $M$ be a compact, oriented 4-manifold.

(i) If $\RRR^2$ is strongly stochastically positive
 and $\pi_1(M)$ is finite, then $M$ is diffeomorphic to $S^4.$

(ii) If $\RRR^2$ is strongly stochastically positive and $\pi_1(M)$ is
infinite, then $\chi (M)\le 0$. In particular, a compact 4-manifold
cannot admit a metric of negative curvature and a metric with $\RRR^2$
strongly stochastically positive.

(iii) If $\RRR^2$ is strongly stochastically positive and $\pi_1(M)$
is nonamenable, then $\chi (M)\le -2$.
\end{proposition}

{\sc Proof:} (i) By \cite{er-inv}, $\RRR^2$
strongly stochastically positive implies $\widetilde M$ is a homotopy
sphere
and hence is diffeomorphic to $S^4$ by Freedman's solution to the
Poincar\'e conjecture.  Since $\chi(\widetilde M)/\chi(M)$ 
equals the order of the covering map from $\widetilde M$ to $M$, the
covering map is either $2-1$ or $1-1$.  Since $M$ is oriented, only
the $1-1$ case is allowed.

(ii) By the $L^2$-index theorem \cite{atiyah}, $\chi (M)$ equals
$L^2\chi
(\widetilde{M})$, the $L^2$ Euler characteristic of $\widetilde M.$  
Since $\pi_1(M)$ is infinite, $0=L^2{
H}^0(\widetilde{M})\cong L^2H^4(\widetilde{M})$. 
Since
$\widetilde{M}$ has $\RRR^2$ strongly stochastically positive for the
pullback
metric, $L^2H^2(\widetilde{M})=0$, and so $0\ge L^2\chi
(\widetilde{M})$.

By Avez's theorem \cite{avez}, the Euler characteristic of a
negatively curved
4-manifold is positive. Of course, such a manifold has infinite
fundamental group.

(iii) As in the proofs of Proposition 2.1 and Theorem 2.3, the
nonamenability of $\pi_1(M)$ implies the existence of a nonconstant
bounded harmonic function $f$ (and hence the transience of Brownian
motion) and the existence of a bounded harmonic one-form $df.$  By
Theorem 2.2,
$0\neq L^2H^1(\widetilde{M})\cong L^2{
H}^3(\widetilde{M})$, so $\chi (M)=L^2\chi (\widetilde{M})
<0$. Since $H^2(M;{\bf R})=0$, $\chi (M)$ is even, so $\chi (M)\le
-2$.\hfill\rule{1.5mm}{3mm}

\bigskip
$S^4$ with its standard metric shows that $\pi_1(M)$ must be infinite
in (ii).
For an example of (iii) not covered by (ii), let $N$ be a compact
hyperbolic 3-manifold with $H^1(N;{\bf R})=0$; examples of such
manifolds occur in Dehn surgery on the figure eight knot complement in
$S^3$. Set $M=N\times S^1$. Then $\chi (M)=0$ and $\pi_1(M)$ is
nonamenable, so $M$ admits no metric with $\RRR^2$ strongly
stochastically
positive. Note that $H^2(M;{\bf R})=0$ and $\pi_2(M)=0$. Thus neither
(i), (ii) nor the Micallef-Moore result \cite{m-m} ($\RRR^2$ pointwise
positive
implies $\pi_2(M^4)=0$) shows that $M$ admits no such metric.

The proposition above is valid for 6-manifolds (except for the part
using Avez's theorem). The
manifolds $N\times N$ and $N\times S^3$, for $N$ as above, do not
admit metrics with $\RRR^2$ strongly stochastically positive, although
this cannot be
seen from the $L^2$-index theorem alone.

\section{Pinching and positivity of $\RRR^p$}

\def\psum{\mathop{{\sum}_p}}

In this section we will give conditions that guarantee stochastic
spectral positivity of ${\cal R}^p.$
\bigskip

\noindent {\bf Definition:}  We say that a Riemannian manifold $M$ is
{\it
$C$-$p$-pinched} if for each $x\in M$ there exists $A= A(x)>0$ such
that for each orthonormal basis
$\{v^1,\ldots,v^n\}$ of $T_x M$,
\begin{equation}\label{sump}
C  A < \sum_{i=1}^p \sum_{j=p+1}^n K(v^i, v^j) < A.\end{equation}
If $C= C(n,p)$ has explicit dependence on $p$, then we will just say
that $M$ is $C$-pinched.  

\bigskip
Here $K(v^i, v^j)$ is the sectional curvature of the plane spanned by
(the duals of) $v^i, v^j.$  We will denote this sum of sectional
curvatures just by $\mathop{\sum_p}$ if the context is clear.

The purpose of this section is to compute $C = C(n,p)$ such that
$C$-pinching implies $M^n$ has $\RRR^p >0.$
Note that a compact manifold is $0$-$1$-pinched iff it has positive
Ricci
curvature, so we are not interested in pinching theorems for
$1$-forms.  Since $\ast \RRR^p =  \RRR^{n-p}\ast$ for the Hodge star
operator $\ast$, we are similarly uninterested in pinching results for
$(n-1)$-forms.  

\begin{theorem} \label{mainpinch} Assume $p\neq 1,n-1.$  Set
$$C(n,p) = \frac{\frac{1}{2} p(n-p) + \frac{4}{3}{p\choose 2}{n-p\choose 2} }
{1+ \frac{1}{2} p(n-p) + \frac{4}{3}{p\choose 2}{n-p\choose 2} }.$$
If $M$ is
$C(n,p)$-pinched, then $\RRR^p >0.$  In particular,
$C(n,p)$-pinching
implies $H^p(M;{\bf R}) =0$ if $M$ is compact.\end{theorem}

Before we begin the proof, we review some curvature tensor
manipulations, both as motivation and for use in a later example.  
 Let $K(i,j)$ denote the sectional curvature of
the 2-plane spanned by $ v^i,v^j$, and similarly let $K(i+j,k)$
denote the sectional curvature for the 2-plane spanned by $
v^i+v^j,v^k$. We have $R_{ijkl}=\langle R(v^i,v^j)v^k,v^l\rangle$,
where we fix our sign convention by $K(1,2) = -R_{1212} = 1$ 
for an orthonormal
frame on $S^2.$
We set $R_{(i+j)k(s+t)l}=\langle
R({v^i+v^j\over \sqrt{2}},v^k)
{v^s+v^t\over \sqq},v^l\rangle$, for $i\neq j, s\neq t$. From
\begin{eqnarray*}
-K(i,l+j) &\equiv& \langle R(v^i,{v^l+v^j\osq})v^i,
{v^l+v^j\osq}\rangle \\
& =&\frac{1}{2}[R_{ilil}+R_{ijij}+R_{ilij}+R_{ijil}],
\end{eqnarray*}
 for $l\neq j$, we get
\begin{equation}\label{rie-one}
R_{ijil}=\left\{ \begin{array}{ll}
-K(i,j),\ \ j=l &\hfil \\
\frac{1}{2}K(i,l)+\frac{1}{2}K(i,j)-K(i,l+j), &\mbox{ $j\not= l$}.
\end{array} \right. 
\end{equation}
We also know
\begin{equation}\label{rie-two}
2R_{(i-k)j(i-k)l}=R_{ijil}-R_{ijkl}-R_{kjil}+R_{kjkl}.
\end{equation}
Switching $i$ and $j$ in (\ref{rie-two}), subtracting the result from
(\ref{rie-two}) and
using the Bianchi identity, we obtain
\begin{equation}\label{rie-three}
3R_{ijkl}=-2R_{(i-k)j(i-k)l}+2R_{(j-k)i(j-k)l}+R_{ijil}
 +R_{kjkl}-R_{jijl}-R_{kikl}.
\end{equation}
Note that combining (\ref{rie-one}) and (\ref{rie-three})
 gives the well known expression for
$R_{ijkl}$ in terms of sectional curvatures.

We now carry these calculations over to $\RRR^p.$  

\begin{lemma} \label{stuff} Let $J, K$ be multi-indices of length
$p$.  

(i) If $|J\cap
K| < p-2$, then $\langle \RRR^p v^J, v^K\rangle =0.$  

(ii)  If $I$ is a multi-index of length $p-2$ and $i,j,k,l$ are
distinct indices not in $I$, then
\begin{equation}\label{one-notes}
\langle R^p( v^i\wedge v^j\wedge v^I), v^k\wedge
v^l\wedge v^I\rangle = 2R_{ijkl}.\end{equation}  

(iii) If $I$ is a multi-index of length $p-2$, if $i,j,l$ are distinct
indices not in $I$, and if $B<\psum <A$, then
\begin{equation}\label{two-notes}
|\langle \RRR^p(v^i\wedge v^j\wedge v^I), v^i\wedge v^l\wedge
v^I\rangle| \leq \frac{1}{2}(A-B). \end{equation}
\end{lemma}

\noindent{\sc Proof:}  Let $a^i$ denote interior multiplication by
$v^i$, and
let $(a^i)^\ast$ denote the adjoint action of wedging with $v^i.$
Then $\RRR^p = R_{ijkl} (a^i)^\ast a^j(a^k)^\ast a^l$
\cite{witten}.      In particular, 
$\RRR^p v^J = a_Kv^K$ has $a_K =0$ if 
$J,K$ differ by more than two indices.  This proves (i).

For (ii), it is immediate that for fixed $p,q,r,s$ we have
$$\langle R_{pqrs} (a^p)^\ast a^q(a^r)^\ast a^s (v^i\wedge v^j\wedge
v^I), v^k\wedge v^l\wedge
v^I\rangle =0,$$
unless $\{p,r\} = \{k,l\}, \{q,s\} = \{i,j\}.$  Now a direct
calculation of the four possibilities for $p,q,r,s$ gives 
\begin{eqnarray*}
\langle R^p( v^i\wedge v^j\wedge v^I), v^k\wedge
v^l\wedge v^I\rangle &=& R_{ljki} + R_{jkli} - R_{likj} + R_{kilj}\\
&=& 2R_{ijkl},\end{eqnarray*}
by the symmetries of curvature tensor.

For (iii), recall that 
\begin{equation}\label{prim}
 \langle \RRR^p(v^1\wedge\ldots\wedge v^p), v^1\wedge\ldots\wedge
v^p\rangle = \sum_{i=1}^p \sum_{j=p+1}^n K(v^i, v^j),\end{equation}
where $\{v^1,\ldots,v^n\}$ is an orthonormal basis extending 
$\{v^1,\ldots,v^p\}$ \cite{EL-RO94}.
 We have
\begin{eqnarray}\label{3.4-old}\lefteqn{
2\langle\RRR^p(v^i\wedge v^j\wedge v^I), v^i\wedge v^l\wedge
v^I\rangle}\\
 &=& \langle \RRR^p [v^i \wedge {v^j+v^l\osq}\wedge v^I],
v^i \wedge {v^j+v^l\osq}\wedge v^I\rangle\nonumber\\
&& - \langle \RRR^p[v^i \wedge {v^j-v^l\osq}\wedge v^I],
v^i \wedge {v^j-v^l\osq}\wedge v^I\rangle.\nonumber
\end{eqnarray}
From the hypothesis, we have
\begin{equation}\label{bottom of eight}
B<\langle \RRR^p[v^i\wedge {v^j+v^l\osq}\wedge v^I],v^i\wedge
{v^j+v^l\osq}\wedge v^I\rangle<A,
\end{equation}
and similarly for 
$\langle \RRR^p[v^i\wedge{v^j-v^l\osq}\wedge v^I],v^i\wedge
{v^j-v^l\osq}\wedge v^I\rangle.$
Combining (\ref{3.4-old}) and (\ref{bottom of eight}) gives (iii).
\hfill\rule{1.5mm}{3mm}
  
\bigskip
We now begin the proof of the theorem.
By (ii) of the Lemma, we may treat $\langle \RRR^p(v^i\wedge v^j\wedge
v^I), v^k\wedge v^l \wedge v^I\rangle$ much like the curvature
tensor.                           
Fix $p, I$ and set
$$\langle R^p( v^i\wedge v^j\wedge v^I), v^k\wedge
v^l\wedge v^I\rangle = T_{ijkl},\ \ T_{ijij} = -S(i,j),$$
\begin{eqnarray*} \langle R^p(v^i\wedge 
 {v^j+ v^m\osq}\wedge v^I), v^i\wedge 
 {v^j+ v^m\osq} \wedge v^I\rangle &=& -S(i,j+m)\\
&=& T_{i(j+m)i(j+m)},\end{eqnarray*}
for distinct indices $i,j,k,l,m\not\in I.$
Note that  $T_{ijkl}= 2R_{ijkl}$,
but the $S$'s are not sectional curvatures --  in fact, $-S(i,j) =
\mathop{\sum_p}$ for some basis $\{v^k\}.$  

Fix $x\in M$ and assume that we have
$B<\psum<A$
at $x$, so
\begin{equation}\label{onea}
B< -S(i,j+m) = T_{i(j+m)i(j+m)} <A.\end{equation}
Of course, $T$ has the same Bianchi identity and symmetries as the
curvature tensor $R$.
Now the argument in (\ref{rie-one})-(\ref{rie-three}) carries over to
$T, S$. In particular,
from
\begin{eqnarray*}
-S(i,l+j) &=& \langle R^p(v^i\wedge {v^l+v^j\osq}
\wedge v^I,{v^l+v^j\osq}\wedge v^I\rangle \\
& =&\frac{1}{2}(T_{ilil}+T_{ijij}+T_{ilij}+T_{ijil}),
\end{eqnarray*}
we get
\[
T_{ijil}=\left\{ \begin{array}{ll}
-S(i,j),\ \ j=l &\hfil \\
\frac{1}{2}S(i,l)+\frac{1}{2}S(i,j)-S(i,l+j), &\mbox{ $j\not= l$}.
\end{array} \right. 
\]
We also know
\[
2T_{(i-k)j(i-k)l}=T_{ijil}-T_{ijkl}-T_{kjil}+T_{kjkl}.
\]
As with $R_{ijkl}$, we obtain
\begin{equation}\label{four}
3T_{ijkl}= -2T_{(i-k)j(i-k)l}+2T_{(j-k)i(j-k)l}+T_{ijil}
 +T_{kjkl}-T_{jijl}-T_{kikl},
\end{equation}
so (\ref{two-notes}) becomes
\begin{equation}\label{five}
|T_{ijil}| \leq \frac{1}{2}(A-B).\end{equation}
Using the estimates (\ref{onea}), (\ref{five}), we see that
(\ref{four}) implies 
\begin{equation} \label{six}
|T_{ijkl}| \leq \frac{4}{3}(A-B).\end{equation}
In summary, for $J\neq K$, we have by Lemma \ref{stuff}(i), (\ref{five})
(\ref{six}),
\[|\langle \RRR^p v^J, v^K \rangle | \leq 
\left\{ \begin{array}{ll}
\frac{1}{2}(A-B),& |J\cap K| = p-1,  \\
\frac{4}{3}(A-B), & |J\cap K| =p-2, \\
0, & |J\cap K| \leq p-3.
\end{array} \right. \]

Thus
\begin{eqnarray}\label{startproof}
\langle \RRR^p(a_Iv^I),a_Jv^J\rangle & = &\sum_I a_I^2\langle
\RRR^pv^I,v^I\rangle +\sum_{I\not= J}a_Ia_J\langle
\RRR^pv^I,v^J\rangle \\
& > &B\sum_I a_I^2-\frac{1}{2}(A-B)\sum_{|I\cap J| =
p-1}|a_Ia_J|\nonumber\\
&&\qquad+
\ft(A-B)\sum_{|I\cap J|=p-2}|a_Ia_J|. \nonumber \end{eqnarray}
  We now estimate the last two terms.

\begin{lemma} \label{choose}
Let $I,J$ be multi-indices of length $p$. Then for $1\leq k\leq p$,
$$\sum_I \sum_{{J\atop |I\cap J| = k}} |a_Ia_J| \leq {p\choose p-k}{n-p\choose 
p-k}\sum_{I}|a_I|^2.$$
\end{lemma}

\noindent{\bf Remark:}  This estimate is sharp, as can be seen by setting $a_I = {n\choose p}^{-1/2}$ for all $I$.

\bigskip

\noindent {\sc Proof:}  For fixed value of $\sum_I a_I^2,$ the 
maximum of the left hand
side of the inequality is attained at some vector $a =
(a_I)$ with $a_I >0.$  Thus we may drop the absolute value signs on
the left hand side.
Set $N = {n\choose p}$ and consider the
$N\times N$ matrix $A$ given by
\[ A_{IJ} = \left\{ \begin{array}{ll}
1,& |I\cap J| = k,  \\
0, & {\rm otherwise}.
\end{array} \right. \]
Here $A$ acts on ${\bf R}^N$ with coordinates indexed by multi-indices
of length $p$.  Finding the maximum of the left hand side for fixed 
$\sum_I a_I^2$ is equivalent to finding the maximum eigenvalue of
$A$.  Observe that 
$$\lambda_0 = {p\choose n-k}{n-p\choose p-k}$$
is an eigenvalue of $A$ with eigenvector $a$ having $a_I = {n\choose
p}^{-1/2}$ for all $I$.

If $p=k$ the result is trivial.  Otherwise we can apply
Perron-Frobenius theory.  Indeed, $A$ has non-negative entries and is
irreducible.  (For if $Q$ is any subset of multi-indices such that
$I\in Q$ and $J\not\in Q$ implies $A_{IJ} =0$, then clearly $I\in Q$
implies $I^\prime\in Q$ if $I^\prime$ is obtained from $I$ by changing
one index; iterating this argument shows that $Q$ consists of all
multi-indices.)  Perron-Frobenius theory (see e.g.~\cite[Theorem 2,
p. 53, and Remark 3, p. 62]{gantmacher}) assures us that such $A$ has
a unique eigenvalue whose eigenvector has positive entries and that
this is the maximum eigenvalue, as we require.
\hfill\rule{1.5mm}{3mm}
\bigskip

Thus (\ref{startproof}) becomes
$$\langle \RRR^p(a_Iv^I),a_Jv^J\rangle  > 
B\Vert a_Iv^I\Vert^2 -\frac{1}{2}(A-B)p(n-p)\sum_I a_I^2 -
\ft(A-B){p\choose 2}{n-p\choose 2}\sum_I a_I^2,$$
and this last term is positive provided 
$$\left(1+\frac{1}{2}p(n-p) +
\ft{p\choose 2}{n-p\choose 2}  \right)B - \left(\frac{1}{2}p(n-p) +
\ft{p\choose 2}{n-p\choose 2}\right)A \geq
0.$$
This finishes the proof of Theorem \ref{mainpinch}.

\bigskip

\noindent {\bf Example:}  By Theorem \ref{cannot}, a compact $n$-manifold
cannot admit both a metric of $[1-(n-2)^{-2}]$-pinched 
negative curvature and a 
$1 - (2/(2+ 2(n-2) +\frac{4}{3}(n-p)(n-p-1)))$-pinched 
metric.  Similar remarks apply to Proposition
\ref{fourdims}.

\bigskip

\begin{corollary}  If
for all $x\in M$, (\ref{sump}) holds for
$C = C(x), A = A(x)$ with
\begin{equation}  A(x) \left( \left[1 + \frac{1}{2}p(n-p)
+\ft{p\choose 2}{n-p\choose 2}\right]C(x) - \left[\frac{1}{2}p(n-p) +
\frac{4}{3}{p\choose 2}{n-p\choose 2}\right]
\right)\label{ssp-pinch}\end{equation} 
strongly stochastically positive, then $\RRR^p $ is strongly
stochastically positive.
In particular, if $M$ is compact
and (\ref{ssp-pinch}) holds, then $H^p(M;{\bf R}) =0.$ \end{corollary}

\noindent {\sc Proof:}  The proof shows that 
$CA<\sum_{i=1}^p \sum_{j=p+1}^n K(v^i, v^j) <A$
 implies 
$$\RRR^p >CA - \frac{1}{2}p(n-p)(A-CA) -
\ft{p\choose 2}{n-p\choose 2}(A-CA).$$
${}$ \hfill\rule{3mm}{3mm} 

\bigskip

We now construct a manifold $M^6$ which is $C(6,3)$-pinched,
and so has $H^3(M;{\bf R})=0$, but which has planes of
negative sectional curvature at each point. 
In contrast,
the Gallot-Meyer
result \cite{g-m} for the curvature operator ${\cal C},$
$${\cal C}>0\ \Rightarrow\ R^k >0,
$$
implies positive sectional curvature, as 
$$\langle {\cal
C}(v^i\land v^j),\  v^i\land v^j\rangle =
\langle -R_{ijkl} v^k\land v^l,\  v^i\land v^j\rangle
=-2R_{ijij}
$$
for an orthonormal frame $\{  v^i\}$. 
Also, we know of no examples
of manifolds with positive curvature operator on isotropic complex
two-planes
(and so $\pi_k(M)=0$ for $1<k\le n/2$ by \cite{m-m})
 but allowing planes of negative
sectional curvature.

We set $M=\Sigma_a\times S^4$, where $\Sigma$ is a closed surface of
genus $g>1$ with a metric of constant negative curvature $-a$, and
$S^4$ is the 4-sphere of constant positive curvature
$1$. We give $M$ the product metric. Since $H^k(M;{\bf
R})\not= 0$ for $k\not= 3$, $M$ cannot be $C(6,k)$-pinched for $k\not=
3$; this is easily verified by computing $\sum_{i=1}^k \sum_{j=k+1}^n
K(v^i,v^j)$ for various permutations of an orthonormal frame $\{
v^i\}$ with $\{ v^1,v^2\}\in T^\ast\Sigma,\ \{ v^3,v^4,v^5,v^6\}\in
T^\ast S^4$.

We give two proofs that $M$ has $\RRR^3 >0$, one proof based on
Theorem \ref{mainpinch}, and one by a direct calculation.

First, for the product metric on $M$, the curvature two-forms
$\Omega_{ij}=R_{ijkl}v^k\land v^l$ vanish  for $\{v^i\}$ as above, 
unless $i,j\in\{ 1,2\}$ or
$i,j\in \{ 3,4,5,6\}$. Moreover, $R_{ijkl}=0$ unless $i,j,k,l\in\{
1,2\},\ \mbox{\rm or}\ i,j,k,l\in\{ 3,4,5,6\}$. Now $R_{1212} =a$, and
for $s,t,r\in\{ 3,4,5,6\}$, (\ref{rie-one}) gives
$$R_{stsr}=\left\{ \begin{array}{ll}
0, &\mbox{$t\not= r$} \\
-1, &\mbox{$t =r$}.
\end{array}\right.
$$
Finally, if $i,j,k,l\in\{ 3,4,5,6\}$ are distinct indices, then
(\ref{rie-three})
implies $R_{ijkl}=0$. This determines the curvature tensor on $M$.

In particular, $R_{ijkl} =0$ if there are at least three distinct
indices among $i,j,k,l.$  From (\ref{one-notes}), 
we see that $\langle R^p v^J,
v^K\rangle =0$ if $|J\cap K| \leq p-2$.  A similar computation shows that
if $|J\cap K| = p-1$, then $\langle R^p v^J,
v^K\rangle$ is a sum of ${R_{ijkl}}$ terms, each of which has three
distinct indices. Thus in this case, $\langle R^p v^J,
v^K\rangle =0$ also.  In summary, we have $\langle \RRR^3 v^I,
v^J\rangle =0$ unless $I=J.$

Let $\{w^i\}$ be an orthonormal frame on $M$.
Writing $w^i = a^i_jv^j$ for an orthogonal
matrix $(a^i_j)$, we get by (\ref{prim})
\begin{eqnarray}\label{twelve}\lefteqn{\sum_{i=1}^3
\sum_{j=4}^6 K(w^i,w^j)}\\ &=&
\sum_{k_1,k_2,k_3,s_1,s_2,s_3}a^1_{k_1}a^2_{k_2}a^3_{k_3} 
a^1_{s_1}a^2_{s_2}a^3_{s_3}\langle \RRR^3(v^{k_1}\wedge v^{k_2}\wedge
v^{k_3}),v^{s_1}\wedge v^{s_2}\wedge
v^{s_3}\rangle\nonumber\\
&=& \sum_{k_1,k_2,k_3}\left(  (a^1_{k_1}a^2_{k_2}a^3_{k_3})^2
\langle \RRR^3(v^{k_1}\wedge v^{k_2}\wedge
v^{k_3}), v^{k_1}\wedge v^{k_2}\wedge
v^{k_3}\rangle\right).\nonumber  \end{eqnarray}
Now 
\begin{equation}\label{thirteen-notes}
\langle \RRR^3(v^{k_1}\wedge v^{k_2}\wedge
v^{k_3}), v^{k_1}\wedge v^{k_2}\wedge
v^{k_3}\rangle = \sum_{i=1}^3
\sum_{j=4}^6 K(\tilde v^i,\tilde v^j),\end{equation}
for some rearrangement $\{\tilde v^i\}$ of $\{ v^i\}$.
Computing the cases where one, two or three of the first three
of the $\tilde v^i$ are in
$TS^4$, we see that the right hand side of (\ref{thirteen-notes}) equals
either $4-a$ or $3$.  Thus if we fix $a=1$, the right hand side of
(\ref{twelve}) becomes 
\begin{eqnarray*}\sum_{i=1}^3
\sum_{j=4}^6 K(w^i,w^j) &=&
 3 \sum_{k_1,k_2,k_3}
(a^1_{k_1}a^2_{k_2}a^3_{k_3})^2\\
 &=& 3\left(\sum_{k_1}(a^1_{k_1})^2\right)
\left(\sum_{k_2}(a^2_{k_2})^2\right)\left(\sum_{k_3}(a^3_{k_3})^2\right)
= 3.\end{eqnarray*}
Thus for $a=1$, $M$ has $\mathop{\sum_3}$ constant, so by the theorem
 $\RRR^3 >0.$  

For the second proof, we compute $R^3$ directly, using our calculation
that $\langle \RRR^3
v^I, v^J\rangle =0$ unless $I=J.$
Let $I = (i_1, i_2, i_3)$, and choose $i_4, i_5, i_6$ so
that $\{i_q: q = 1,...,6\} = \{1,...,6\}$.
  By (3.7),
we see that 
\begin{equation}\label{ten-notes}
\langle \RRR^3(a_Iv^I), a_Iv^I\rangle = (\sum_I a_I^2) \langle \RRR^3
v^I, v^I\rangle = 
\sum_{I, |I|=3} a_I^2 \sum_{q=1}^3 \sum_{s=4}^6
K(v^{i_q},v^{i_s}).\end{equation}
Since $\mathop{\sum_3}$ equals either $4-a$ or $3$, taking 
$a<4$ makes the right hand side of (\ref{ten-notes}) positive.  

\bigskip
\noindent{\bf Remark:} Let $\{w^i\}$ be an orthonormal frame.
We can use [8,Thm. 5A] to conclude
the stronger result $H_3(M;{\bf Z})=0$ once we check 
\begin{equation}\label{eleven-notes} \sum_{i=1}^3
\sum_{j=4}^6 K(w^i,w^j) >\frac{\Vert\alpha\Vert^2}{2}
-\frac{n|H|^2}{2} \end{equation}
pointwise on $M$, where $\alpha$ is the second
fundamental form, $H$ is the mean curvature for an isometric
immersion of $M$ in some ${\bf R}^N$, and $n=6$. Here we first isometrically
immerse $\Sigma_a$ in some ${\bf R}^{N_1}$ for a fixed $a$ and $S^4$
in ${\bf R}^5$ in the usual way. We then put $M$ isometrically in
${\bf R}^{N_1+5}$.

For $a\leq 1$, we know that the left hand side of (\ref{eleven-notes})
is at least $3$.
On the other hand, for $a$ small,
$\Vert\alpha\Vert^2$ is very close to $\Vert\alpha\Vert^2$ for $S^4$,
and similarly for $n|H|^2$. On $S^4$, $\Vert\alpha\Vert^2=4$
and $H$ is minus the radial vector, so on $M$,
 $(\Vert\alpha\Vert^2 /2)-(n|H|^2/2)$ is approximately
$2/3$.
 Thus (\ref{eleven-notes}) is satisfied.


\renewcommand{\thefootnote}{}

\def\limsup{\mathop{\overline{\rm lim}}}
\def\liminf{\mathop{\underline{\rm lim}}}
\def\exp{{\rm e}}

\section{Remarks on spectral gap estimates}

The crucial  Lemma 2.1 treats the $L^\infty$ theory of harmonic
forms, so it is natural to look for an
$L^2$ version.
For this, let $\lambda_1^q=\lambda_1^{h,q}$ be the spectral gap at zero 
for $q$-forms, i.e. 
$\lambda_1^q=\inf \left\{ \hbox{spec}(\Delta^{h,q})-\{0\}\right\}$.
 Equivalently
$$\lambda_1^q=\inf_{\phi\in C_K^\infty, \phi\in H_q^{\perp}}
 {\int_M \langle\Delta^h \phi,\phi\rangle \, \mu(dx)
\over |\phi|_{L^2}^2}.$$
Here $H_q^\perp$ is  the space of differential $q$-forms perpendicular
 to the the harmonic forms.  We will show that the $L^2$ version of Lemma
 2.1 corresponds to  lower bound estimates  for $\lambda_1^q$
 for differential forms.

\begin{proposition}
 Let $1\le q\le n-1$, and suppose $\R^{q\pm1}$ are bounded from
 below. Then
$$\lambda_1^q \ge \min\left\{\lambda_1^{q-1}, \lambda_1^{q+1}\right\}.$$
\end{proposition}

\noindent{\sc Proof:}  
For  $\phi, \psi\in C_c^\infty \Lambda^q$, and in the notation 
of \S2, we have

\begin{eqnarray*}
&&\left |\, \int_M \langle  P_t\phi-\phi,\psi \rangle \, \mu(dx) \, \right |\\
&=&\half \left| \int_0^t \int_M \langle  d \phi, P_s(d\psi)\rangle 
 \mu(dx) \, ds
+ \int_0^t \int_M  \langle  \delta^h\phi,
 P_s(\delta^h\psi)\rangle  \mu(dx) \, ds	\right|\\
&=& \half \left|  \int_M \langle  d \phi,\int_0^t P_s(d\psi) \, ds\,\rangle  
\mu(dx) \right|
+\half \left| \int_M \langle  \delta^h\phi,
\int_0^t P_s(\delta^h\psi) \, ds\,\rangle  \mu(dx)	\right|\\
&\le& \half ||d\phi||_{L^2}\hskip 2pt \cdot \hskip 2pt
  ||  \int_0^t P_s(d\psi) ds||_{L^2}+
 \half ||\delta^h \phi||_{L^2} \hskip 2pt \cdot \hskip 2pt
 ||\int_0^t P_s(\delta^h \psi) ds||_{L^2}\\
 \end{eqnarray*}

Suppose $\lambda_1^{q+1}>0$ and $\lambda_1^{q-1}>0$. Then since $d\psi$
and $\delta^h \psi$ are orthogonal to the $L^2$ h-harmonic forms, we have 
the existence of $\displaystyle{\int_0^\infty P_s(d\psi) ds}$
and  $\displaystyle{\int_0^\infty P_s(\delta^h \psi) ds}$ in $L^2$ with

$$||\int_0^\infty P_s(d\psi) ds||\le
{2\over \lambda_1^{q+1}} \, \cdot\, ||d\psi||_{L^2}$$

and
$$||\int_0^\infty P_s(\delta^h \psi) ds|| \le
{2 \over \lambda_1^{q-1} }\, \cdot \,  ||\delta^h \psi||_{L^2}.
$$

Thus letting $t\to \infty$ in the earlier inequality gives

$$\left|\int_M\langle H\phi-\phi,\psi\rangle \mu(dx)\right |
\le {1\over \lambda_1^{q+1}} \left(|d \phi|_{L^2} \right)\, \cdot \,
 \left(|d \psi|_{L^2} \right)
+ {1\over \lambda_1^{q-1}} \left( |\delta^h \phi|_{L^2} \right)
\, \cdot \,\left(|\delta^h \psi|_{L^2} \right).$$

Taking $\phi=\psi\in H^\perp$ and using $H\phi=0$, we have

$$\left( ||\phi||_{L^2}\right)^2 \le
  { 1\over \hbox{min}\left\{\lambda_1^{q-1}, \lambda_1^{q+1}\right\}}
\,\cdot\, \langle\Delta \phi,\phi\rangle,$$
and so
 $\lambda_1^q \ge \hbox{min}\left\{\lambda_1^{q-1}, \lambda_1^{q+1}\right\}.$
\hfill\rule{3mm}{3mm}

\bigskip

If $M$ is compact, this result is an easy consequence of Hodge theory,
as $d$ and $\delta$ give a conjugacy
of the restriction of $\Delta^q$ to $({\rm Ker}\ 
 d)^\perp $ with $\Delta^{q+1}$
restricted to the image of $d$, and similarly for $\Delta^q$ on
 $ ({\rm Ker}\ \delta)^\perp$.  However, this argument does not extend
to the noncompact case without more stringent
curvature conditions, due to the possible presence of
continuous spectrum.

We can also estimate $\lambda_1^q$ in terms of $\underline{ \RRR}^{q\pm
1}$, but these estimates
are useful only when there are no $L^2$ harmonic $(q\pm 1)$-forms.
If we denote by  $\lambda_0(f)$  the bottom of the
 spectrum of the operator  $\Delta+ f$ on functions, then we have
\[
\lambda_1^q\ge \min\left\{ \lambda_0(\underline{R}^{q-1}),
\lambda_0(\underline{R}^{q+1})\right\},
\]
assuming  $\RRR^{q+1}$ and $\RRR^{q-1}$  are bounded from below.
This comes from \cite[Theorem 3C]{EL-RO88} 
 which states that the bottom
of the spectrum of $\Delta^p$ is bigger or equal to
 $\lambda_0(\underline{R}^{p})$ for each $p\in \{1,2,\dots, n-1\}$.
Note that by \cite[Proposition 4B]{EL-RO88}, we have
$$\lambda_0(\underline{\RRR}^p) \geq -2 \overline{\lim_{t\to\infty}}
\frac{1}{t} \log {\Bbb E} e^{-\frac{1}{2}\int_o^t
\underline{\RRR}^p(x_s)ds},$$ 
which gives 
 an estimate for $\lambda_1^q$ in terms of $\underline{\cal R}^{q-1}$
and   $\underline{\cal R}^{q+1}$.

\bibliographystyle{amsplain}
\bibliography{paper5}

\ifx\undefined\bysame
\newcommand{\bysame}{\leavevmode\hbox to3em{\hrulefill}\,}
\fi
\begin{thebibliography}{10}

\bibitem{atiyah}
M.~F. Atiyah, {\em Elliptic operators, discrete groups, and von {N}eumann
  algebras}, Asterisque {\bf 32-33} (1976), 43--72.

\bibitem{avez}
A.~Avez, {\em Applications de la formule de {G}auss-{B}onnet-{C}hern aux
  vari\'et\'es \`a quatre dimensions}, C. R. Acad. Sc. Paris {\bf 256} (1963),
  5488--5490.

\bibitem{azencott}
R.~Azencott, {\em Behavior of diffusion semigroups at infinity}, Bull. Math.
  Soc. France {\bf 102} (1974), 193--240.

\bibitem{bakry}
D.~Bakry, {\em Un crit\'ere de non-explosion pour certaines diffusions sur une
  vari\'et\'e {R}iemannienne compl\'ete}, C. R. Acad. Sc. Paris {\bf 303}
  (1986), 23--26.

\bibitem{brooks}
R.~Brooks, {\em The fundamental group and the spectrum of the {L}aplacian},
  Comment. Math. Helv. {\bf 56} (1981), 581--598.

\bibitem{el-li94}
K.~D. Elworthy and X.-{M}. Li, {\em Formulae for the derivatives of heat
  semigroups}, J. Functional Analysis {\bf 125} (1994), 252--286.

\bibitem{EL-RO88}
K.~D. Elworthy and S.~Rosenberg, {\em Generalized {B}ochner theorems and the
  spectrum of complete manifolds}, Acta Appl. Math. {\bf 12} (1988), 1--33.

\bibitem{er-inv}
\bysame, {\em Manifolds with wells of negative {R}icci curvature}, Inventiones
  Math. {\bf 103} (1991), 471--495.

\bibitem{EL-RO93}
\bysame, {\em The {W}itten {L}aplacian on negatively curved simply connected
  manifolds}, Tokyo J. Math. {\bf 16} (1993), no.~2, 513--524.

\bibitem{EL-RO94}
\bysame, {\em Homotopy and homology vanishing theorems and the stability of
  stochastic flows}, Geometric and Functional Analysis {\bf 6} (1996), 51--78.

\bibitem{gaffney}
M.~P. Gaffney, {\em The conservation property of the heat equation on
  {R}iemannian manifolds}, Comm. Pure Appl. Math. {\bf 12} (1959), 1--11.

\bibitem{g-m}
S.~Gallot and D.~Meyer, {\em Op\'erateur de courbure et {L}aplacien des formes
  diff\'erentielles d'une vair\'et\'e {R}iemannienne}, J. Math. Pures et Appl.
  {\bf 54} (1975), 259--284.

\bibitem{gantmacher}
F.~R. Gantmacher, {\em The {T}heory of {M}atrices}, vol.~2, Chelsea Publishing
  Co., New York, 1960.

\bibitem{lawson-simons}
H.~B. Lawson, Jr. and J.~Simons, {\em On stable currents and their applications
  to global problems in real and complex geometry}, Annals of {M}ath. {\bf 98}
  (1973), 427--450.

\bibitem{Li.thesis}
X.-{M}. Li, {\em Stochastic flows on noncompact manifolds}, Ph.D. thesis,
  University of Warwick, 1992.

\bibitem{group}
\bysame, {\em On extension of {M}yers' theorem}, Bull. London Math. Soc. {\bf
  27} (1995), 392--396.

\bibitem{ls}
T.~Lyons and D.~Sullivan, {\em Function theory, random paths and covering
  spaces}, J. {D}ifferential {G}eometry {\bf 19} (1984), 299--323.

\bibitem{mckean}
H.~McKean, {\em An upper bound for the spectrum of ${\Delta}$ on a manifold of
  negative curvature}, J. Differential Geometry {\bf 4} (1970), 359--366.

\bibitem{m-m}
M.~Micallef and J.~D. Moore, {\em Minimal two-spheres and the topology of
  manifolds with positive curvature on totally isotropic two-planes}, Annals of
  Math. {\bf 87} (1988), 199--227.

\bibitem{m-w}
M.~Micallef and J.~Wolfson, {\em The second variation of area of minimal
  surfaces in four manifolds}, preprint.

\bibitem{reed-simon}
M.~Reed and B.~Simon, {\em Methods of {M}odern {M}athematical {P}hysics},
  Academic Press, Orlando, 1980.

\bibitem{ry}
S.~Rosenberg and D.~Yang, {\em On the fundamental group of manifolds of almost
  positive {R}icci curvature}, J. Math. Soc. Japan {\bf 46} (1994), 267--288.

\bibitem{witten}
E.~Witten, {\em Supersymmetry and {M}orse theory}, J. Differential Geometry
  {\bf 121} (1985), 169--186.

\bibitem{yau}
S.-T. Yau, {\em Harmonic functions on complete {R}iemannian manifolds}, Comm.
  Pure and Applied Math. {\bf 38} (1975), 201--228.

\end{thebibliography}

\end{document}